\documentclass[submission,copyright,creativecommons]{eptcs}
\providecommand{\event}{AiML 2026} 

\usepackage{graphicx}

\usepackage{mathUtil}

\usepackage{aiml26}

\usepackage{iftex}

\ifpdf
  \usepackage{underscore}         
  \usepackage[T1]{fontenc}        
\else
  \usepackage{breakurl}           
\fi

\title{Uniform Local Tabularity in Intuitionistic Logic}
\author{Rodrigo Nicolau Almeida
\institute{ILLC - UvA - Amsterdam - Netherlands}
\email{r.dacruzsilvapinadealmeida@uva.nl}
}

\newcommand{\titlerunning}{Uniform Local Tabularity in Intuitionistic Logic}
\newcommand{\authorrunning}{R.N.~Almeida}

\hypersetup{
  bookmarksnumbered,
  pdftitle    = {\titlerunning},
  pdfauthor   = {\authorrunning},
  pdfsubject  = {EPTCS},               
  pdfkeywords = {keyword1, keyword2} 
}

\date{}

\begin{document}

\maketitle

\begin{abstract}
By contrast with  $\mathsf{S4}$, the analysis of local tabularity above $\mathsf{IPC}$ has provided a difficult challenge. This paper studies a strengthening of local tabularity --- \textit{uniform local tabularity} --- where one demands that all formulas be equivalent to formulas of a given implication depth. Algebraically, this amounts to considering Heyting algebras generated by finitely many iterations of the implication operation. It is shown that in contrast with locally finite Heyting algebras, $n$-uniformly locally finite Heyting algebras always form a variety, and an explicit axiomatization of the variety of $n$-uniform locally finite Heyting algebras for $n\leq 3$ is given. In connection with this analysis, it is shown that there exist locally tabular logics which are not uniformly locally tabular, answering a question of Shehtman ---- an example of a pre-uniformly locally tabular logic is presented, which is shown to be the unique pre-uniformly locally tabular extension of the system $\mathsf{KG}$.
\end{abstract}

\section{Introduction}

In intuitionistic and modal logics, one often looks for properties which ensure the system is ``sufficiently simple'' to analyze in depth. Local tabularity ensures this by saying that over finitely many propositional variables there can only be finitely many such formulas, and in turn, that the semantics of the system cannot be too complicated. This appears as a strengthening of `tabularity' --- saying that a single finite structure captures the full semantics of the logic --- which is quite well-understood in intuitionistic and modal logics (see \cite[Chapter 12]{Chagrov1997-cr}).  Moreover, it is a classical result due to Segerberg \cite{Segerberg1971-SEGAEI} and Maksimova \cite{Maksimova1975} that an extension of $\mathsf{S4}$ is locally tabular if and only if it is of finite depth. This provides an easy criterion and algorithm to verify local tabularity, with there being a single ``pre-locally tabular'' logic \cite[Sec.~12.4]{Chagrov1997-cr}.

By contrast, in the system $\mathsf{IPC}$, the situation is quite different. It was shown by Mardaev \cite{Mardaev1984,Mardaev1987} that there is a continuum of pre-locally tabular logics, all above $\mathsf{KC}$. It is open whether every locally tabular logic is contained in a pre-locally tabular one and whether local tabularity is decidable. In terms of algebraic models --- which correspond in algebraic logical terms to the logics --- it is thus open whether the category of locally finite Heyting algebra admits any easy and logically insightful description (see \cite{HYTTINEN2024} for recent discussion).

Given these difficulties it is reasonable to seek for strengthenings of local tabularity. In \cite{shehtmanrieger}, Shehtman observed that one can use bisimulation games to characterize when certain logics are locally tabular, introducing a property called ``finite formula depth'', which corresponds to the existence of a bound on the nestings of modal operators or implications. In \cite{colimitsofHeytingalgebras}, we dubbed this property ``uniform local tabularity'', with the idea being that the bound on modal depth or implication depth provides a uniform bound on the size of formulas. It was shown by Shehtman that above $\mathsf{S4}$ local tabularity implies uniform local tabularity, via finite depth, and hence the `niceness' of local tabularity may be attributed to the uniform version. The relationship between local tabularity and its uniform version was likewise noted in \cite{shapirovskyfiniteheight}, where generalizations of the finite height criterion were discussed. The concept of uniform local tabularity also has deep connections to the theory of free Heyting algebras and free constructions of Heyting algebras: as discussed in \cite[Chapter 5]{almeidaphdthesis}, uniformity of local tabularity appears to be connected to the existence of simple representations of duals of free algebras.

This article provides a first systematic investigation of this property in the context of superintuitionistic logics. Our main contributions are the following:
\begin{enumerate}
    \item We provide an explicit axiomatization of the varieties of $n$-uniform Heyting algebras for $n\leq 3$; for $n=1$, this corresponds to the logic $\mathsf{CPC}$, and for $n=2$, the extensions of logic $\mathsf{LC}$, the G\"{o}del-Dummett calculus \cite{Dummett1959,davis1990}; for $n=3$, such a logic is introduced here, and called  $\mathsf{3Uni}$, being the least $3$-uniform logic, and sitting between well-known logics:
    \begin{equation*}
        \mathsf{wPL}\subsetneq \mathsf{3Uni}\subsetneq \mathsf{LC},
    \end{equation*}
    where $\mathsf{wPL}$ is the logic of the Weak Peirce's Law  (see Section \ref{Section: 1-uniform and 2-uniform logics}).
    \item We provide an example of a locally tabular logic which is not uniformly locally tabular, answering a problem left open in \cite{Shehtman2016}. Indeed, we show that the logic $\mathsf{Box}$ (see Section \ref{sec: preuniform local tabularity}), is pre-uniformly locally tabular --- all of its extensions are uniformly locally tabular, though $\mathsf{Box}$ itself is not so. We also show that it is the unique logic with this property above the system $\mathsf{KG}$ \cite{kuznetsovgerciu}.
\end{enumerate}

\section{Preliminaries and basic theory}

\subsection{Superintuitionistic logics and Kripke models}

Throughout we will assume that we are working in the signature $\mathcal{L}=(\wedge,\vee,\rightarrow,\bot,\top)$ of Heyting algebras, the language of intuitionistic logic. We use the letters $\phi, \psi,\dots$ for $\mathcal{L}$-formulas and letters $p, q,\dots$  or $x,y, \dots$ for variables; we compactly represent a tuple of distinct variables
as $\overline{p}$. We write $\varphi(\overline{p})$ to mean that all the variables of $\varphi$ occur amongst $\overline{p}$.

\begin{definition}\label{def: literal implication depth}
    The \textit{implication-depth} of $\varphi$, denoted by $d(\varphi)$, is defined by induction:
\begin{itemize}
    \item $d(p)=d(\bot)=d(\top)=0$;
    \item $d(\varphi\wedge\psi)=d(\varphi\vee\psi)=\max(d(\varphi),d(\psi))$;
    \item $d(\varphi\rightarrow\psi)=\max(d(\varphi),d(\psi))+1$.
\end{itemize}
\end{definition}

Given a poset $P$, we write $P\Vdash \varphi$ to mean that $P$ validates $\varphi$ in the usual Kripke semantics (see e.g., \cite[Chapter 3]{Chagrov1997-cr}). We denote by $\mathrm{Ext}(\mathsf{IPC})$ the lattice of its extensions; elements are called \textit{superintuitionistic logics}. Given a poset $P$, we denote by:
\begin{equation*}
    \mathrm{Log}(P)=\{\varphi : P\Vdash \phi\},
\end{equation*}
the \textit{logic of $P$}, which is an element of $\mathrm{Ext}(\mathsf{IPC})$. If $P\Vdash L$ for all $\varphi\in L$ we say that $P$ \textit{is a frame for $L$}.

Given a poset $P$, a subset $U\subseteq P$ is called an \textit{upset} if whenever $x\in U$ and $x\leq y$ then $y\in U$. We denote the set of upsets of $P$ by $\mathsf{Up}(P)$. Given a subset $C$, we denote by ${\uparrow}C$ the least upset containing it, and for $x\in P$, ${\uparrow}x\coloneqq {\uparrow}\{x\}$. We say that $Q$ is a \textit{rooted upset} of $P$ if there is $x\in P$ such that $Q= {\uparrow}x$. Given $x\in P$, we also denote ${\uparrow}x$ by $P_{x}$.

\begin{definition}
    A \textit{p-morphism} between posets $f:P\to Q$ is an order-preserving function such that for each $x\in P$ and $y\in Q$ if $f(x)\leq y$ then there is some $x'\geq x$ such that $f(x')=y$; we say that $Q$ is a \textit{p-morphic image} of $P$ if there is a surjective p-morphism from $P$ onto $Q$.
\end{definition}

We now recall some relevant notions concerning finite Kripke models; the reader can find all necessary background in \cite[Chapter 4.1]{Ghilardi2002}.

\begin{definition}
    We say that a triple $(\bbM,x)=(P,\leq,c,x)$ where $(P,\leq)$ is a rooted poset, with root $x$, and $c:P\to \mathcal{P}(\overline{p})$ is an order-preserving function, is a \textit{model over the propositional letters $\overline{p}$}. We say that $(\bbM,x)$ is a model of a logic $L$ (over the propositional letters $\overline{p}$) if $(P,\leq)$ is a frame for $L$. Moreover, given $\mathbf{K}$ a class of posets, we say that $(\bbM,x)$ a model over $\overline{p}$ is a \textit{$\overline{p}$-model of $\mathbf{K}$} if there is some poset $Q\in \mathbf{K}$ such that $P$ is a p-morphic image of a rooted upset of $Q$. We denote by $\mathcal{M}_{\overline{p}}(\mathbf{K})$ the class of $\overline{p}$-generated models over $\mathbf{K}$. We also denote by $M_{\mathbf{K}}(\varphi)=\{(\mathfrak{M},x)\in \mathcal{M}_{\overline{p}}(\mathbf{K}) : (\mathfrak{M},x)\Vdash\varphi\}$.
\end{definition}

\begin{definition}
    Given two models $(\mathfrak{M},x)=(P,\leq,c,x)$ and $(\mathfrak{N},y)=(Q,\leq,c',y)$ over $\overline{p}$, we say that they are \textit{bisimilar} if there exists a relation $S\subseteq P\times Q$ such that (1) $(x,y)\in S$; (2) Whenever $wSw'$ then $c(w)=c'(w')$; (3) whenever $wSw'$ and $w\leq v$ there is some $w'\leq v'$ such that $vSv'$; (4) whenever $wSw'$ and $w'\leq v'$ there is some $w\leq v'$ such that $vSv'$. We say that two points $w\in P$ and $w'\in Q$ are bisimilar if $(P_{w},w)$ is bisimilar with $(Q_{w'},w')$.
\end{definition}

\begin{remark}
    \ 
    \begin{enumerate}
        \item (\textbf{Notation}) Whenever we do not wish to emphasize the specific propositional letters used in the language,we will  identify the set $\{p_{1},\dots,p_{n}\}$ with the natural number $n$, and consider models over $\overline{p}$ to be maps $c:P\to 2^{n}$; the image for a point $x\in P$ will be called its \textit{color}, denoted by a tuple $c(x)=\langle i_{1},\dots,i_{n}\rangle$. We will further shorten tuples of the form $\langle 1,\dots,1,0,\dots,0\rangle$, consisting of a sequence of $1$ followed by a sequence of $0$, by the notation $1^{n-k}0^{k}$. When $n$ is given, we will often omit the order and the valuation, and write simply $(P,x)$ for a model. We occasionally write $\mathcal{M}_{n}(\mathbf{K})$ for the class of models over $n$ from $\mathbf{K}$. 
        \item All statements given below for posets $(P,x)$ or models thereof assume that $P$ is \emph{rooted at $x$}. This is implicit in our notion of model, but we repeat it for the sake of the reader familiar with frame theoretic arguments.
        \item Given an $n$, say that a model $(P,x)$ over $n$ is \textit{n-generated} if $P$ does not contain distinct bisimilar points. A well-known characterization of $n$-generation which we will use throughout is through $\alpha$ and $\beta$-reductions (see \cite[Chapter 3]{bezhanishviliphdthesis}). Namely, given a model $P$, we say that we can apply an \emph{alpha-reduction} ($\alpha$) whenever $x\in P$ has an immediate successor $x'$ such that ${\uparrow}x=({\uparrow}x'\cup\{x\})$ and $c(x)=c(x')$, obtained by identifying $x$ with $x'$; we say we can perform a $\beta$-reduction ($\beta$) if $x,x'$ are two points, such that ${\uparrow}x-\{x\}={\uparrow}x'-\{x'\}$, and $c(x)=c(x')$, obtained by identifying $x$ and $x'$. A model is $n$-generated precisely when there are no $\alpha$ or $\beta$-reductions which can be performed on the model. A convenient fact which we will exploit throughout is that for $n$-generated models, \emph{bisimulation implies isomorphism}.
    \end{enumerate}
\end{remark}

\begin{definition}
    Let $(P,\leq,c, x)$ and $(Q,\leq,c',y)$ be two models over $n$. We say that a relation $S_{k}\subseteq P\times Q$ is a \textit{k-bisimulation} based on $(x,y)$ if there are relations $S_{k}\subseteq \dots\subseteq S_{j}\subseteq\dots\subseteq S_{0}\subseteq P\times Q$ for each $0\leq j\leq k$, such that $xS_{k}y$ and for each $w\in P,w'\in Q$:
    \begin{enumerate}
        \item Whenever $wS_{0}w'$ then $c(w)=c'(w')$;
        \item For each $j<k$, whenever $wS_{j+1}w'$ and $w\leq v$ there is some $w'\leq v'$ such that $vS_{j}v'$.
        \item For each $j<k$, whenever $wS_{j+1}w'$ and $w'\leq v'$ there is some $w\leq v$ such that $vS_{j}v'$.
\end{enumerate}
We say that $(P,x)$ is $k$-bisimilar with $(Q,y)$ if there exists a $k$-bisimulation between them. We write $(P,x)\leq_{0} (Q,y)$ if $c(x)\leq c(y)$ and $(P,x)\leq_{k+1}(Q,y)$ if for each $y\leq w'$ there is some $x\leq w$ such that $(P_{w},w)\sim_{k}(Q_{w'},w')$. We also write $(P,x)\leq_{\infty}(Q,y)$ if whenever $y\leq w'$ there is some $x\leq w$ such that $(P,w)$ is bisimilar with $(W,w')$
\end{definition}

As described in \cite[Chapter 4.1]{Ghilardi2002}, $n$-bisimulations can be characterized in game-theoretic terms: given a pair $(P,x)$ and $(Q,y)$, Player I wins a $0$-bisimulation game if $c(x)\neq c(y)$ and otherwise Player II wins; and Player I wins an $n+1$-bisimulation game if there is $z\in P$, such that Player II cannot choose any $w\in Q$ and win an $n$-bisimulation game with $(P_{z},z)$ and $(Q_{w},w)$.

The following is a basic property of $n$-bisimulation:

\begin{proposition}\label{prop: invariance of truth under upwards}
    For each pair of models $(P,x)$ and $(Q,y)$ over $\overline{p}$ the following are equivalent:
    \begin{enumerate}
        \item $(P,x)\leq_{n}(Q,y)$;
        \item Whenever $d(\varphi(\overline{p}))\leq n$ and $(P,x)\Vdash \varphi$ then $(Q,y)\Vdash\varphi$.
    \end{enumerate}
\end{proposition}

It follows from this that for finite models, $(P,x)\leq_{\infty}(Q,y)$ if and only if for every formula $\phi$, $(P,x)\Vdash \phi$ then $(Q,y)\Vdash \phi$. The following is completely analogous to \cite[Chapter 4, Exercise 11]{Ghilardi2002}, where we relativise to a specific subclass of posets, to capture different equivalence in different logics:

\begin{proposition}\label{prop: Equivalence condition for uniform local finiteness}
    Let $\mathbf{K}$ be a class of finite posets, $L=\mathrm{Log}(\mathbf{K})$ and let $\overline{p}$ be a set of propositional variables. Then we have:
    \begin{enumerate}
        \item\label{eq: fmp property} For each $\varphi(\overline{p}),\psi(\overline{p})$, we have:
        \begin{equation*}
            L\vdash \varphi\leftrightarrow\psi \iff \mathcal{M}_{\mathbf{K}}(\varphi)=\mathcal{M}_{\mathbf{K}}(\psi).
        \end{equation*}
        \item \label{eq: equivalence of formulas} The following are equivalent:
    \begin{enumerate}
        \item $S\subseteq \mathcal{M}_{\overline{p}}(\mathbf{K})$ is the class of models from $\mathcal{M}_{\overline{p}}(\mathbf{K})$ satisfying a formula $\psi$ with $d(\psi)\leq n$;
        \item $S$ is closed under $\leq_{n}$-bisimulation within $\mathcal{M}_{\overline{p}}(\mathbf{K})$: whenever $(P,x)\in S$ and and $(Q,y)\in \mathcal{M}_{\overline{p}}(\mathbf{K})$ and $(P,x)\leq_{n}(Q,y)$ then $(Q,y)\in S$ 
    \end{enumerate}
    \end{enumerate}
\end{proposition}
\begin{proof}
    \eqref{eq: fmp property} is an immediate consequence of $L$ being the logic of a class of finite posets, and therefore, if $\phi\leftrightarrow\psi\notin L$, there being a model of $\phi$ which is not a model of $\psi$.

    \eqref{eq: equivalence of formulas} If $S$ is the class of models satisfying a formula $\psi$ with $d(\psi)\leq n$, then by Proposition \ref{prop: invariance of truth under upwards}, we have that $S$ is closed under $\leq_{n}$-bisimulation. Conversely, assume that $S$ is closed under $\leq_{n}$-bisimulation, and consider:
    \begin{equation*}
        \mathsf{Th}_{n}(S)=\{\psi(\overline{p}) : d(\psi)\leq n, \forall (P,x)\in S, (P,x)\Vdash \psi\}.
    \end{equation*}
    Since there are only finitely many formulas of implication depth at most $n$, $\mathsf{Th}(S)$ is finite, and so we can consider $\chi(S)=\bigwedge\mathsf{Th}(S)$. We claim that $S$ is precisely the class of models from $\mathcal{M}_{\overline{p}}(\mathbf{K})$ satisfying $\chi(S)$. It is clear that if $(P,x)\in S$, $(P,x)\Vdash \chi(S)$. On the other hand, assume that $(Q,y)\notin S$. Then for each $(P,x)\in S$, $(P,x)\nleq_{n}(Q,y)$, which by Proposition \ref{prop: invariance of truth under upwards} means that there is a formula $\phi(\overline{p})$ of implication depth at most $n$ such that $(P,x)\Vdash \phi_{P}$ and $(Q,y)\nVdash \phi_{P}$. Let $\mu_{P}=\bigvee_{P\in S}\phi_{P}$; note that up to equivalence, $\mu_{P}$ is finite, and $\mu_{P}\in \mathsf{Th}_{n}(S)$. Since $(Q,y)\nVdash \mu_{P}$, then $(Q,y)\nVdash \chi(S)$, which concludes the proof.
\end{proof}

\subsection{Local tabularity, local finiteness and their uniform versions}

To investigate these notions it will be convenient to work with Heyting algebras. For the connection between superintuitionistic logics and varieties of Heyting algebras, and the algebraic semantics therein, we refer to \cite[Chapter 8]{Chagrov1997-cr}. Moreover, though the reader will not explicitly need any facts about the duality of finite Heyting algebras and finite posets, this will be implicit in the following results.

Given a formula $\phi$, we write $H\vDash \phi$ to mean that the Heyting algebra $H$ validates $\phi$; notations are extended to classes in the usual way. Given a class of Heyting algebras $\mathbf{K}$, we write:
\begin{equation*}
    \mathrm{Log}(\mathbf{K})=\{\phi : \mathbf{K}\vDash \phi\};
\end{equation*}
conversely, given a logic $L$, we write:
\begin{equation*}
    \mathrm{Alg}(L)=\{H\in \mathbf{HA} : H\vDash L\}.
\end{equation*}

We now turn to the basic concepts in our analysis.

\begin{definition}
    Let $L\supseteq \mathsf{IPC}$ be a superintuitionistic logic. We say that $L$ is:
    \begin{itemize}
        \item \textit{locally tabular} if given $\overline{p}=\{p_{1},\dots,p_{n}\}$ propositional variables, there are only finitely many formulas over $\overline{p}$  up to $L$-provability. We denote by $\mathsf{LTab}$ the class of locally tabular logics;
        \item \textit{n-uniform} if each formula $\phi$ is equivalent to a formula of implication depth $n$ over $L$. We say that $L$ is \textit{uniformly locally tabular} if $L$ is $n$-uniform for some $n$. We denote by $\mathsf{ULTab}_{n}$ the class of $n$-uniform logics, and by $\mathsf{ULTab}$ their union, the class of all uniformly locally tabular logics;
        \item of \textit{finite depth} if $L\supseteq\mathsf{BD}_{n}$ for some $n$, where:

        \begin{equation*}
            \mathsf{bd}_{0}=p_0 \text{ and } \mathsf{bd}_{n+1}=p_{n+1}\vee(p_{n+1}\to\mathsf{bd}_{n}),
        \end{equation*} and   \begin{equation*}
            \mathsf{BD}_{n}=\mathsf{IPC}\oplus \mathsf{bd}_{n}.
        \end{equation*}
        We denote by $\mathsf{FinDep}$ the class of logics of finite depth.
        \item \textit{tabular} if there is a finite poset $P$ such that $L=\mathrm{Log}(P)$. We write $\mathsf{Tab}$ for the class of tabular logics.
    \end{itemize}
\end{definition}

All of these concepts correspond algebraically to specific properties:

\begin{definition}
    Let $H$ be an algebra. We say that $H$ is \textit{locally finite} if whenever $A\leq H$ is a finitely generated algebra, then $A$ is finite. We say that a class of algebras $\mathbf{K}$ is \textit{locally finite} if every algebra in $\mathbf{K}$ is locally finite.
\end{definition}

Recall that given $H$ a Heyting algebra, $A\subseteq H$ a subset, the least Heyting subalgebra of $H$ containing $A$ can be constructed as:
\begin{equation*}
    \langle A\rangle =\{t(b_{1},\dots,b_{k}) : b_{i}=a_{i} \text{ or } b_{i}=\neg a_{i}, a_{i}\in A, \ t(x_{1},\dots,x_{n}) \text{ is a formula }\}.
\end{equation*}
Local finiteness means precisely that for each $H$, and $A$ a finite subset, there is an $n$ such that $\langle A\rangle$ is generated by formulas of implication depth at most $n$. This is where we can impose uniformity:

\begin{definition}
    Let $H$ be an algebra. We say that $H$ is $n$-\textit{uniformly finite}, or simply \textit{n-uniform}, if whenever $A\subseteq H$ is finite, then:
    \begin{equation*}
        \langle A\rangle=\{t(a_{1},\dots,a_{k}) : a_{i}\in A, t(x_{1},\dots,x_{k}) \text{ has implication depth $\leq n$}\}.
    \end{equation*}
    In this case we say that $\langle A\rangle$ is \textit{n-uniformly generated}. We say that a class of algebras $\mathbf{K}$ is $n$-\textit{uniformly} finite if every algebra in $\mathbf{K}$ is $n$-uniformly finite. 
\end{definition}

The motivation for the former definition lies in a certain modal understanding of Heyting algebras: one sees the implications as modalities and their nesting as a way to parametrize the complexity of a given logical system or algebra. A different way of measuring such complexity goes through the structure of frames:

\begin{definition}
    Given $x\in P$ a poset, the \textit{depth of $x$}, denoted $d_{p}(x)$, is the largest $n$ such that there exists a chain $x=x_{1}\leq \dots\leq x_{n}$, and there is no such chain of size $n+1$, or $\infty$ if such an $n$ does not exist. If each point in $P$ has a depth, we denote by $d_{p}(P)=\sup\{d_{p}(x) : x\in P\}$.
\end{definition}

The following is all well-known and provides the algebraic translation of the above properties:

\begin{proposition}\label{prop: Inclusion of classes of logics}
    The following holds for $L$ a superintuitionistic logic:
    \begin{enumerate}
        \item $L$ is locally tabular if and only if $\mathsf{Alg}(L)$ is locally finite.
        \item $L$ is $n$-uniform if and only if $\mathsf{Alg}(L)$ is $n$-uniformly finite.
        \item $L$ is of finite depth $k$ if and only if whenever $P\Vdash L$ is a poset, then $P$ is of depth at most $k$, if and only if $L$ contains the formula $\mathsf{bd}_{k}$, defined recursively by:
        \begin{equation*}
            \mathsf{bd}_{1}=p_{1}\vee \neg p_{1} \text{ and } \mathsf{bd}_{j+1}=p_{j+1}\vee (p_{j+1}\rightarrow \mathsf{bd}_{j})
        \end{equation*}
        \item $L$ is tabular if and only if $\mathsf{Alg}(L)$ is finitely generated as a variety (i.e., it is the least variety containing a given finite algebra).
    \end{enumerate}
    Moreover, the following inclusions hold between classes of Heyting algebras:
    \begin{equation*}
        \mathsf{Tab}\subseteq \mathsf{FinDep}\subseteq \mathsf{ULTab}\subseteq \mathsf{LTab}.
    \end{equation*}
\end{proposition}
\begin{proof}
    The equivalences of (1-3) are easy to prove. We concentrate on the inclusions. $\mathsf{Tab}\subseteq \mathsf{FDep}$, since each finite poset $P$ satisfies $\mathsf{bd}_{k}$ where $k$ is its height; the fact that $\mathsf{FDep}\subseteq \mathsf{ULTab}$ was shown by Shehtman \cite{Shehtman2016}, and amounts to the fact that for posets of finite depth, there is a uniform bound for how complex $n$-bisimulation can be; the fact that $\mathsf{ULTab}\subseteq \mathsf{LTab}$ is obvious.
\end{proof}

Most converses of the inclusions of Proposition \ref{prop: Inclusion of classes of logics} are known not to hold. Namely, the logic $\mathsf{BD}_{2}$ is known not to be tabular, even though it is of finite depth. The logic $\mathsf{LC}$ is axiomatized as:
\begin{equation*}
    \mathsf{LC}\coloneqq (p\rightarrow q)\vee (q\rightarrow p),
\end{equation*}
is known to be uniformly locally tabular (see \cite{Shehtman2016}, see also \cite{colimitsofHeytingalgebras} for a proof) but is not of finite depth. We will take a look at the last converse --- $\mathsf{LTab}\subseteq \mathsf{ULTab}$ in Section \ref{section: hierarchy of locally tabular logics}.

\subsection{Basic properties and criteria for uniform local tabularity}

In order to establish $n$-uniform local tabularity we will make use of the following characterization in terms of $n$-bisimulations\footnote{We note for the reader that the criterion below makes use of the relation $\leq_{n}$ as opposed to $n$-bisimilarity. In \cite{Shehtman2016}, the corresponding result does not mention this relation, however, we have not been able to reproduce such a result when omitting this relation.}:

\begin{proposition}\label{Saturated classes yielding n-uniformity}
    Let $\mathbf{K}$ be a class of finite posets. Then $L=\mathsf{Log}(\mathbf{K})$ is $n+1$-uniformly locally tabular if and only if for each $\overline{p}$, whenever $(P,x)$ and $(Q,y)$ are two models in $\mathcal{M}_{\overline{p}}(\mathbf{K})$ and $(P,x)\leq_{n+1}(Q,y)$, then $(P,x)\leq_{\infty}(Q,y)$. Moreover, if $n$-bisimilarity implies full bisimilarity, the latter holds.
\end{proposition}
\begin{proof}

If $L$ is $n+1$-uniformly locally tabular, then all formulas are equivalent to formulas of implication depth $\leq n+1$. Let $(P,x),(Q,y)\in \mathcal{M}_{\overline{p}}(\mathbf{K})$ be two models, and assume that $(P,x)\leq_{n+1}(Q,y)$. Then whenever $(P,x)\Vdash \phi$ for any $\phi$, $(Q,y)\Vdash \phi$. Thus $(P,x)\leq_{\infty}(Q,y)$ by our remark after Proposition \ref{prop: invariance of truth under upwards}.

Conversely, let $\varphi(\overline{p})$ be an arbitrary formula. Consider $S=\{(P,x)\in \mathcal{M}_{\overline{p}}(\mathbf{K}) : (P,x)\Vdash \phi\}$. Now assume that $(P,x)\leq_{n+1}(Q,y)$. By assumption, $(P,x)\leq_{\infty}(Q,y)$, and consequently, $(Q,y)\Vdash \phi$ as well. This means that means that $S$ is closed under $\leq_{n+1}$-bisimulation within $\mathcal{M}_{\overline{p}}(\mathbf{K})$, which by Proposition \ref{prop: Equivalence condition for uniform local finiteness}.(2) that $S$ is the class of models from $\mathcal{M}_{\overline{p}}(\mathbf{K})$ satisfying a formula $\psi$ with $d(\psi)\leq n+1$. By Proposition \ref{prop: Equivalence condition for uniform local finiteness}.(1), this means that $L\vdash \varphi\leftrightarrow\psi$, i.e., $\varphi$ is equivalent to a formula of implication depth $\leq n+1$. Thus $L$ is $n+1$-uniformly locally tabular, as desired.

If we assume that $n$-bisimilarity implies full bisimilarity instead, then if we assume that $(P,x)\leq_{n+1}(Q,y)$, by considering $y\leq y$, we have that there is some $x\leq w$ such that $(P_{w},w)\sim_{n}(Q,y)$. Since $(P,x)\Vdash \phi$, then $(P_{w},w)\Vdash \phi$ as well; since $n$-bisimilarity implies full bisimilarity, by invariance of truth under bisimilarity, $(Q,y)\Vdash \phi$. We thus obtain the same conclusion. 
\end{proof}

We close this section by noting an interesting property of $n$-uniform Heyting algebras. Let $\mathbf{HA}_{lf}$ be the category of locally finite Heyting algebras with Heyting algebra homomorphisms, and let $\mathbf{HA}_{lf}^{n}$ be the category of $n$-uniformly finite Heyting algebras with the same maps.

\begin{proposition}\label{prop: closure of n-uniformity under variety}
    The category $\mathbf{HA}_{lf}$ is closed under subalgebras, homomorphic images and finite products. For each $n\in \omega$, $\mathbf{HA}_{lf}^{n}$ is a variety. 
\end{proposition}
\begin{proof}
    The fact about subalgebras is trivial. To see the fact about homomorphic images, say $H$  is locally finite (n-uniformly finite), and $A$ is a homomorphic image. Let $\{b_{1},\dots,b_{n}\}\subseteq A$ be finitely many elements. Then we can find finitely many elements $\{h_{1},\dots,h_{n}\}\subseteq H$, mapping surjectively onto $b_{1}$; by restricting the homomorphism from $H$ to $A$, we find that $\langle \{h_{1},\dots,h_{n}\}\rangle$ maps surjectively to $\langle \{b_{1}\dots,b_{n}\}\rangle$. Thus the latter is finite ($n$-uniformly generated) by hypothesis.

    To see the fact about products, note that if $H_{0},H_{1}$ are locally finite, and $(a_{0},b_{0}),\dots,(a_{n},b_{n})$ where $a_{i}\in H_{0}$ and $b_{i}\in H_{1}$, then
    \begin{equation*}
        \langle (a_{0},b_{0}),\dots,(a_{n},b_{n})\rangle = \{t((a_{0},b_{0}),\dots,(a_{n},b_{n})) : d(t) \text{ is the maximum of that of $\langle A\rangle$ and $\langle B\rangle$.}
    \end{equation*}
    The same argument works for infinite products if we assume the algebras are $n$-uniformly finite, since then, rather than taking the maximum of $\langle A\rangle$ and $\langle B\rangle$ one can take the bound $n$.
\end{proof}

The above Proposition entails that there is the least $n$-uniform logic for each $n$. On the other hand, $\mathbf{HA}_{lf}$ is not a variety ---  it contains every finite Heyting algebra, and hence the smallest variety containing it is $\mathbf{HA}$ itself.

\section{The hierarchy of locally tabular logics}\label{section: hierarchy of locally tabular logics}

In this section we will study the hierarchy of locally tabular logics --- the \textit{top} of the Kuznetsov hierarchy \cite{Bezhanishvili2019}. In \ref{Section: 1-uniform and 2-uniform logics} by characterizing the $0,1,2,3$-uniformly locally tabular logics. We then give in Section \ref{section: non-uniform locally tabular logic} an example of a locally tabular but not uniformly locally tabular logic.

\subsection{0-uniform, 1-uniform, 2-uniform, 3-uniform logics}\label{Section: 1-uniform and 2-uniform logics}

The reader will notice that due to our definition, the only superintuitionistic logic which is $0$-uniform is the trivial logic. $1$-uniformity turns out to be quite a restrictive property: only $\mathsf{CPC}$ has it. $2$-uniformity will be characteristic of extensions of $\mathsf{LC}$, whilst $3$-uniformity and beyond will turn out to have a rather complex structure.

We assume the reader is familiar with the Yankov theorem (see \cite[Chapter 3]{bezhanishviliphdthesis} for an extensive discussion) associating a formula $\mathcal{J}(Q)$ to a rooted poset $Q$; the key property we will need is the following:

\begin{theorem}[Yankov's theorem]
    For a finite poset $P$ and a rooted poset $Q$, $P\nVdash \mathcal{J}(Q)$ if and only if $Q$ is a rooted subframe of a p-morphic image of $P$.
\end{theorem}

We briefly recall the structure of the lattice of extensions of superintuitionistic logics; all proofs of these facts can be found in \cite[Chapter 4, Chapter 8]{Chagrov1997-cr}:
\begin{enumerate}
    \item Given logics $L'\subseteq L$ we say that $L'$ is a \emph{cover} of $L$ if $L$ is an immediate successor of $L'$, i.e., $L'\subseteq K\subseteq L$ implies that either $L'=K$ or $K=L$.
    \item $\mathsf{CPC}$ is the largest consistent extension of $\mathsf{IPC}$;
    \item $\mathsf{CPC}$ has a unique cover $\mathsf{Sme}=\mathsf{IPC}\oplus (\neg p\vee q\vee (q\rightarrow p))=\mathsf{Log}(2)$, the logic of the $2$-element poset.
    \item $\mathsf{Sme}$ is an extension of $\mathsf{LC}$, in fact the greatest extension different from $\mathsf{CPC}$. 
    \item The logic $\mathsf{LC}$ is axiomatized by $\mathcal{J}(\mathsf{2F})$ and $\mathcal{J}(\mathsf{Dia})$ where $\mathsf{2F}$ is the 2-fork (the frame underlying the model $M_{2}$ in Figure \ref{fig:PnandQnposets}), and $\mathsf{Dia}$ is the diamond frame (the poset $Q_{6}$ in Figure \ref{fig:qiposetsfor2uniformity}).
\end{enumerate}

Then we have:

\begin{theorem}
    For $L$ a superintuitionistic logic, we have:
    \begin{enumerate}
        \item $L$ is $0$-uniform if and only if $L$ is inconsistent.
        \item $L$ is $1$-uniform if and only if $L=\mathsf{CPC}$;
        \item $L$ is $2$-uniform if and only if $L\supseteq \mathsf{LC}$.
    \end{enumerate}
\end{theorem}
\begin{proof}
As noted, the only $0$-uniform logic is the inconsistent one.  Note that $\mathsf{CPC}$ is $1$-uniform: by Definition \ref{def: literal implication depth}, every formula in $\mathsf{CPC}$ is equivalent to a conjunction of disjunctions of formulas $p,\neg p$. On the other hand, the logic $\mathsf{Sme}$ is not $1$-uniform: we see in Figure \ref{fig:PnandQnposets} models $M_{0}$ and $M_{1}$, based on frames validating $\mathsf{Sme}$, and $M_{1}\leq_{1}M_{0}$, certainly, but $M_{1}\nleq_{\infty}M_{0}$, since no point in $M_{1}$ labelled with $0$ is $1$-bisimilar with the singleton containing $0$. By Proposition \ref{Saturated classes yielding n-uniformity}, this means that the logic is not $1$-uniform, and neither is any logic extended by it.

As noted before, it was shown by Shehtman \cite{Shehtman2016} that $\mathsf{LC}$ is $2$-uniform; it is also not difficult to see: if $(P,x)$ and $(Q,y)$ are two linearly ordered posets rooted at $x,y$, they are $n$-generated, and they are $1$-bisimilar, then notice that by $n$-generation, there cannot be successive nodes with the same color (otherwise we could identify them). But then by an easy induction we see that the two chains ${\uparrow}x$ and ${\uparrow}y$ must be isomorphic. It follows that if $(P,x)\sim_{1}(Q,y)$ for $n$-generated models, then $(P,x)$ and $(Q,y)$ are isomorphic, and consequently, fully bisimilar.  We can then appeal to the same proposition as above.

Now assume that $L$ is some logic and $\mathsf{LC}\nsubseteq L$; then $L\nVdash \mathcal{J}(\mathsf{2F})$ or $L\nVdash \mathcal{J}(\mathsf{Dia})$, so either the 2-fork or the diamond are frames for $L$. If $\mathsf{2F}$ is a frame for $L$, note that in Figure \ref{fig:PnandQnposets}, we have the models $M_{1}$ and $M_{2}$, and $M_{2}\leq_{2}M_{1}$: if we pick in $M_{1}$ the top node labelled with $1$, there is a node in $M_{2}$ with the same label, obtaining an isomorphism; and if we pick in $M_{1}$ the root, node that $M_{1}\sim_{1}M_{2}$. However, $M_{2}\nleq_{3}M_{1}$, since by picking the root of $M_{1}$, there is no $2$-bisimilar node in $M_{2}$.   Similarly, the models $M_{3}$ and $M_{4}$ (the frame of which is a p-morphic image of $M_{3}$) have the property that $M_{3}\leq_{2}M_{4}$ (there is an isomorphic node in $M_{3}$ for each node but the root, and the roots are $1$-bisimilar), but $M_{2}\nleq_{3}M_{4}$. This shows that $L$ cannot be $2$-uniformly locally tabular again by Proposition \ref{Saturated classes yielding n-uniformity}.
\end{proof}

\begin{figure}
    \centering
\begin{tikzpicture}
\node at (-2,0) {$\bullet$};
\node at (-2,-0.3) {$0$};
\node at (-2,-1) {$M_{0}$};

    \node at (0,0) {$\bullet$};
    \node at (0,-1) {$M_{1}$};
    \node at (0,-0.3) {$0$};
    \node at (0,1) {$\bullet$};
    \node at (0,1.3) {$1$};

    \node at (3,0) {$\bullet$};
    \node at (3,-1) {$M_{2}$};
    \node at (3,-0.3) {$0$};
    \node at (2,1.3) {$1$};
    \node at (2,1) {$\bullet$};
    \node at (4,1) {$\bullet$};
    \node at (4,1.3) {$0$};

    \draw (0,0) -- (0,1);
    \draw (2,1) -- (3,0) -- (4,1);
\end{tikzpicture}
\qquad
\begin{tikzpicture}
    \node at (0,0) {$\bullet$};
    \node at (0,-0.3) {$00$};
    \node at (1,1) {$\bullet$};
    \node at (1.3,1) {$00$};
    \node at (-1,1) {$\bullet$};
    \node at (-1.3,1) {$10$};
    \node at (0,2) {$\bullet$};
    \node at (0,2.3) {$11$};
    \node at (0,-1) {$M_{3}$};
    \draw (0,0) -- (-1,1) -- (0,2) -- (1,1) -- (0,0);

    \node at (3,0) {$\bullet$};
    \node at (3.3,0) {$00$};
    \node at (3,1) {$\bullet$};
    \node at (3.3,1) {$10$};
    \node at (3,2) {$\bullet$};
    \node at (3.3,2) {$11$};
    \node at (3,-1) {$M_{4}$};

    \draw (3,0) -- (3,1) -- (3,2);
\end{tikzpicture}

\caption{Models witnessing degrees of uniformity}
    \label{fig:PnandQnposets}
\end{figure}

The case of $3$-uniform logics is a fair bit more complicated and deserves further consideration. For that purpose, recall the logic $\mathsf{wPL}$:
\begin{equation*}
    \mathsf{wPL}\coloneqq \mathsf{IPC}\oplus (q\rightarrow p)\vee ((p\rightarrow q)\rightarrow p)\rightarrow p.
\end{equation*}

It is known (see \cite[A.9]{Esakiach2019HeyAlg}) that a Heyting algebra $H$ validates $\mathsf{wPL}$ if and only if it is a \textit{cascade Heyting algebra}. Their variety was named by Esakia $\mathsf{CHL}$, and was studied in \cite{Bezhanishvili2021} in the context of analyzing profiniteness and representability of Heyting algebras. The finite members of such a variety are the ``sums'' of Boolean algebras; we will make use of the dual posets of such algebras:

\begin{definition}
    A finite (rooted) poset $(P,\leq)$ is called a \textit{Boolean sum} if $d_{p}(P)=n$, and:
    \begin{enumerate}
        \item Whenever $x\prec y$, i.e., $y$ is an immediate successor of $x$, then $d_{p}(x)=d_{p}(y)+1$;
        \item Each point $x$ of depth $k+1$ sees all points of depth $k$ for each $k$.
    \end{enumerate}
\end{definition}

When working with Boolean sums, note that the notion of depth of a point captures precisely its \emph{level}: for $P$ a Boolean sum, $\mathsf{Lev}_{k}(P)=\{x\in P : d_{p}(x)=k\}$ consists of incomparable points which are all seen by the points in $\mathsf{Lev}_{k-1}(P)$, and all see the points in $\mathsf{Lev}_{k+1}(P)$. The following is a consequence of \cite[A.9.2]{Esakiach2019HeyAlg}

\begin{proposition}
    A finite rooted poset $P$ is a frame for $\mathsf{wPL}$ if and only it is a Boolean sum.
\end{proposition}

The following was shown by Mardaev \cite[Theorem 5]{Mardaev1991}:

\begin{theorem}\label{Local tabularity of wPL}
    The logic $\mathsf{wPL}$ can be axiomatized as:
\begin{equation*}
    \mathsf{wPL}\coloneqq \mathsf{IPC}\oplus \mathcal{J}(Q_{1})\oplus \mathcal{J}(Q_{2})\oplus \mathcal{J}(Q_{3}),
\end{equation*}
where $\mathcal{J}(Q_{i})$ is the Yankov formula for the posets in Figure \ref{fig:qiposetsfor2uniformity}. Moreover, this logic is locally tabular.    
\end{theorem}

We will be concerned with special extensions of this logic. For that we will need some definitions:

\begin{definition}\label{def: stack depth}
    Let $P$ be a Boolean sum of depth $n$. For each $n$, let $c_{n}$ be the cardinality of the set of elements of depth $n$. Given $\langle c_{1},\dots,c_{n}\rangle$ we say that a consecutive subsequence $c_{i},\dots,c_{i+k}$ where $c_{i+j}>1$ for $0\leq j\leq k$ is a $k$-\textit{stack}. We say that $P$ has \textit{stack-depth} $k$ if there exists a $k$-stack but not a $k+1$-stack in $P$.
\end{definition}

The stack depth of a Boolean sum thus measures how many consecutive levels there can be which contain more than one point.

\begin{remark}
    In \cite{Bezhanishvili2021} a subvariety of $\mathsf{CHL}$ was considered in connection to the problem of profinite completions, called the variety of \textit{Diamond Heyting algebras} $\mathsf{DHL}$. Such algebras are of width $2$, have a maximal element, and they satisfy the Yankov formula $\mathcal{J}(Q_{5})$. Thus, their finite rooted posets are precisely Boolean sums of stack depth $1$, width $2$ having a greatest element (see also Section \ref{sec: preuniform local tabularity}, where the associated logic will be called $\mathsf{Box}$).
\end{remark}

\begin{lemma}\label{lem: Axiomatization Lemma}
Let $(P,x)$ be a finite poset which is a Boolean sum. Then $P\Vdash \mathcal{J}(Q_{i})$, for $i\in \{4,5\}$ in Figure \ref{fig:qiposetsfor2uniformity}, if and only if $P$ has stack-depth $1$.
\end{lemma}
\begin{proof}
First note that both of the posets $Q_{4},Q_{5}$ have stack depth $2$. If $Q_{4},Q_{5}$ is a rooted subframe of $Q$, which is a p-morphic image of $P$, by pulling along the p-morphism we will see that $P$ has stack depth at least $2$.

Conversely, assume that $P$ has stack depth at least $2$. Let $x\in P$ be a maximal point such that ${\uparrow}x$ has stack-depth $2$. We may thus assume that $x$ has at least two successors $y_{0},y_{1}$, and that each of them has two successors (otherwise we could have picked $y_{i}$ as our maximal point) $z_{0},z_{1}$. Then there are two cases: either $z_{0}$ and $z_{1}$ have no common successor, or they do. In the first first case, we define a map to $Q_{4}$ as follows: we send $z_{0}$ to the left top point in $Q_{4}$, $z_{1}$ and all remaining points of that depth to the right top point; $y_{0}$ to the left middle point, $y_{1}$ and all remaining points of that depth to the middle right point; and $x$ to itself. This is a surjective p-morphism onto $Q_{4}$. The case of $Q_{5}$ is similar, except we send all points above $z_{0},z_{1}$ to the maximal point.
\end{proof}

\begin{figure}
    \centering
\begin{tikzpicture}
    \node at (1,0) {$\bullet$};
    \node at (0,1) {$\bullet$};
    \node at (0,2) {$\bullet$};
    \node at (2,1) {$\bullet$};
    \node at (1,-0.5) {$Q_{1}$};

    \draw (0,2) -- (0,1) -- (1,0) -- (2,1);

    \node at (3,1) {$\bullet$};
    \node at (3,2) {$\bullet$};
    \node at (4,0) {$\bullet$};
    \node at (4,-0.5) {$Q_{2}$};
    \node at (5,1) {$\bullet$};
    \node at (5,2) {$\bullet$};

    \draw (3,2) -- (3,1) -- (5,2) -- (5,1) -- (4,0) -- (3,1);

    \node at (6,1.5) {$\bullet$};
    \node at (7,0) {$\bullet$};
    \node at (7,-0.5) {$Q_{3}$};
    \node at (8,1) {$\bullet$};
    \node at (8,2) {$\bullet$};
    \node at (7,3) {$\bullet$};

    \draw (7,3) -- (6,1.5) -- (7,0) -- (8,1) -- (8,2) -- (7,3);

    \node at (4,-1) {$\bullet$};
    \node at (3,-2) {$\bullet$};
    \node at (5,-2) {$\bullet$};
    \node at (4,-3) {$\bullet$};

    \draw (4,-1) -- (3,-2) -- (4,-3) -- (5,-2) -- (4,-1);

    \node at (4,-3.5) {$Q_{6}$};

    \node at (6,-1) {$\bullet$};
    \node at (6,-2) {$\bullet$};
    \node at (7,-3) {$\bullet$};
    \node at (8,-2) {$\bullet$};
    \node at (8,-1) {$\bullet$};
    \node at (7,-3.5) {$Q_{7}$};

    \draw (6,-1) -- (6,-2) -- (7,-3) -- (8,-2) -- (8,-1);

    \node at (9,-2) {$\bullet$};
    \node at (10,-2) {$\bullet$};
    \node at (10,-3) {$\bullet$};
    \node at (11,-2) {$\bullet$};

    \draw (9,-2) -- (10,-3) -- (10,-2) -- (10,-3) -- (11,-2);

    \node at (10,-3.5) {$Q_{8}$};

    \node at (9,1) {$\bullet$};
    \node at (9,2) {$\bullet$};
    \node at (10,0) {$\bullet$};
    \node at (10,-0.5) {$Q_{4}$};
    \node at (11,1) {$\bullet$};
    \node at (11,2) {$\bullet$};

    \draw (11,2) -- (11,1) -- (9,2) -- (9,1) -- (11,2) -- (11,1) -- (10,0) -- (9,1);

    \node at (12,1) {$\bullet$};
    \node at (12,2) {$\bullet$};
    \node at (13,0) {$\bullet$};
    \node at (13,-0.5) {$Q_{5}$};
    \node at (14,1) {$\bullet$};
    \node at (14,2) {$\bullet$};
    \node at (13,3) {$\bullet$};

    \draw (13,3) -- (12,2) -- (13,3) -- (14,2) -- (14,1) -- (12,2) -- (12,1) -- (14,2) -- (14,1) -- (13,0) -- (12,1);
    
\end{tikzpicture}    \caption{$Q_{i}$ posets}
    \label{fig:qiposetsfor2uniformity}
\end{figure}

Given a rooted frame $(P,x)$ based on a Boolean sum with stack depth $1$, let us say that a point $y\in P$ is \emph{solitary} if $y$ is the only point of depth $k$ for some $k$. Otherwise we say that $y$ is \emph{accompanied} by a point $y'$. The importance of stack depth for us comes from the following result:

\begin{lemma}\label{lem: 1-bisim}
Let $(P,x)$ be an $n$-generated model based on a finite Boolean sum of stack depth $1$. If $(P,x)\sim_{1}(P_{y},y)$ for some $y\in {\uparrow}x$, then $x=y$.
\end{lemma}
\begin{proof}
Assume first that $x$ is a solitary point. Then for each $y$, either $x\leq y$ or $y\leq x$ must hold. If $y\leq x$, note that by $n$-generation, any predecessor of $x$ must have a color strictly smaller than $x$, and thus we cannot have $(P,y)\sim_{1}(P,x)$. On the other hand, if $x\leq y$, let $S(x)$ be the set of immediate successors of $x$. If $S(x)$ is a singleton, we likewise have that the color of $x$ must be strictly smaller, and so $(P,y)\sim_{1}(P,x)$ cannot hold. If $S(x)=\{x_{0},x_{0}'\}$, then since the stack depth is $1$, these points must have distinct color, and either $c(x_{0})\neq y$ or $c(x_{0}')\neq y$; say $c(x_{0})\neq y$ without loss of generality. If $(P,y)\sim_{1}(P,x)$, then $y$ would have to have a successor $y\leq y_{0}$ such that $c(y_{0})=c(x_{0})$.  However, note that if $z$ is the unique immediate successor of $x_{0},x_{0}'$, $c(x_{0}),c(x_{0}')<c(z)$, so $y_{0}=x_{0}$ or $y_{0}=x_{0}'$. Thus if $x\neq y$, we would have a contradiction, since $x\leq y\leq y_{0}$ would force $y=y_{0}$ (since $y_{0}$ must be an immediate successor of $x$), whilst we chose $c(y)\neq c(y_{0})$. It follows that $(P,x)\sim_{1}(P,y)$ can only hold if $x=y$.

Next assume that $x$ is accompanied by $x'$. If $x<y$, then by $n$-generation, $c(x)\neq c(y)$. Moreover, $c(x)\neq c(x')$, so $(P,x)\sim_{1}(P,x')$ cannot hold. Finally if $y<x$, then by the same arguments as in the previous case, if $(P,x)\sim_{1}(P,y)$, we would obtain a contradiction. We thus obtain the result.
\end{proof}

\begin{proposition}\label{Characterization of 2-uniformity}
Let $P,Q$ be two finite Boolean sum frames of stack depth $1$. Assume that $(P,x)\sim_{2}(Q,y)$ are $n$-generated models which are $2$-bisimilar. Then they are isomorphic.
\end{proposition}
\begin{proof}
Let $S_{2}\subseteq S_{1}\subseteq S_{0}\subseteq P\times Q$ be a $2$-bisimulation. We claim that $S_{1}$ defines an isomorphism, mapping $x$ to $y$. Indeed, assume that $y,y'\in P$ and $(y,z)\in S_{1}$ and $(y',z)\in S_{1}$. Then we have that $(P,y)\sim_{1}(Q,z)\sim_{1}(P,y')$, so  $(P,y)\sim_{1}(P,y')$. By Lemma \ref{lem: 1-bisim}, we have that $y=y'$. Thus $S_{1}$ is injective; a similar argument applied to $Q$ shows that $S_{1}$ is functional, and the bisimulation clauses guarantee that it is a surjective relation. Thus $S_{1}$ is an isomorphism as desired.
\end{proof}

Let $\mathsf{2Uni}$ be the following logic:

\begin{equation*}
    \mathsf{2Uni}\coloneqq \mathsf{IPC}\oplus \{\mathcal{J}(Q_{i}) : i\in \{1,\dots,5\}\} = \mathsf{wPL}\oplus \mathcal{J}(Q_{4})\oplus \mathcal{J}(Q_{5}).
\end{equation*}

Then we have:

\begin{theorem}\label{thm: characterization of 3-uniform}
    The logic $\mathsf{2Uni}$ is the least $3$-uniform logic.
\end{theorem}
\begin{proof}
    First, note that in Figure \ref{fig:2bisimilar models not bisimilar}, we see that all the posets $Q_{i}$ for $1\leq i\leq 5$ admit $2$-bisimilar models (obtained by p-morphic images of $Q_{i}$, denoted $Q_{i}^{l}$ and $Q_{i}^{r}$ whether they occur on the right or left) which violate our criterion. Namely, we have that $Q_{i}^{l}\leq_{3}Q_{i}^{r}$, since for each node except the root in $Q_{i}^{r}$ there is a point isomorphic to it in $Q_{i}^{l}$, and the root of $Q_{i}^{r}$ is $2$-bisimilar to $Q_{i}^{l}$; but $Q_{i}^{l}\nleq_{4}Q_{i}^{r}$, because the root in $Q_{i}^{r}$ has no $3$-bisimilar point on the left models.

    Hence, no logic which contains them can be $2$-uniform. By our axiomatization of the logic $\mathsf{2Uni}$, any logic which is not an extension of it will therefore contain one of these posets by Yankov's theorem.
    
    On the other hand, by Lemma \ref{lem: Axiomatization Lemma}, Proposition \ref{Characterization of 2-uniformity} and Corollary \ref{Saturated classes yielding n-uniformity}, we have that $\mathsf{2Uni}$ is $3$-uniform.
\end{proof}

\begin{figure}
        \centering
        \begin{tikzpicture}

        \node at (-2,-1) {$Q_{1}$};
            \node at (0,0) {$\bullet$};
            \node at (0,0.3) {$1$};
            \node at (0,-1) {$\bullet$};
            \node at (-0.3,-1) {$0$};
            \node at (1,-2) {$\bullet$};
            \node at (1,-2.3) {$0$};
            \node at (2,-1) {$\bullet$};
            \node at (2.3,-1) {$0$};

            \draw (0,0) -- (0,-1) -- (1,-2) -- (2,-1);

            \node at (3,-1) {$\bullet$};
            \node at (3,-0.7) {$1$};
            \node at (4,-2) {$\bullet$};
            \node at (4,-2.3) {$0$};
            \node at (5,-1) {$\bullet$};
            \node at (5,-0.7) {$0$};

            \draw (3,-1) -- (4,-2) -- (5,-1);

            \node at (7,-1) {$Q_{2}$};

            \node at (9,0) {$\bullet$};
            \node at (9,0.3) {$11$};
            \node at (9,-1) {$\bullet$};
            \node at (8.7,-1) {$10$};
            \node at (10,-2) {$\bullet$};
            \node at (10,-2.3) {$00$};
            \node at (11,0) {$\bullet$};
            \node at (11,0.3) {$10$};
            \node  at (11,-1) {$\bullet$};
            \node at (11.3,-1) {$00$};

            \draw (9,0) -- (9,-1) -- (11,0) -- (11,-1) -- (10,-2) -- (9,-1);

            \node at (12,0) {$\bullet$};
            \node at (12,0.3) {$11$};
            \node at (13,-1) {$\bullet$};
            \node at (13,-0.7) {$10$};
            \node at (14,0) {$\bullet$};
            \node at (14,0.3) {$10$};

            \node at (13,-2) {$\bullet$};
            \node at (13,-2.3) {$00$};

            \draw (13,-2) --(13,-1) -- (12,0) -- (13,-1) -- (14,0);

            \node at (-2,-6) {$Q_{3}$};

            \node at (7,-6) {$Q_{4}$};

            \node at (1,-4) {$\bullet$};
            \node at (1,-3.7) {$11$};
            \node at (0,-5.5) {$\bullet$};
            \node at (-0.3,-5.5) {$00$};
            \node at (1,-7) {$\bullet$};
            \node at (1,-7.3) {$00$};
            \node at (2,-6) {$\bullet$};
            \node at (2.3,-6) {$00$};
            \node at (2,-5) {$\bullet$};
            \node at (2.3,-5) {$10$};

            \draw (1,-4) -- (0,-5.5) -- (1,-7) -- (2,-6) -- (2,-5) -- (1,-4);

            \node at (3,-6) {$\bullet$};
            \node at (3.3,-6) {$00$};
            \node at (4,-5) {$\bullet$};
            \node at (4,-4.7) {$11$};
            \node at (4,-7) {$\bullet$};
            \node at (4,-7.3) {$00$};
            \node at (5,-6) {$\bullet$};
            \node at (5.3,-6) {$10$};

            \draw (3,-6) -- (4,-7) -- (5,-6) -- (4,-5) -- (3,-6);

            \node at (9,-5) {$\bullet$};
            \node at (9,-4.7) {$11$};
            \node at (11,-5) {$\bullet$};
            \node at (11,-4.7) {$10$};
            \node at (9,-6) {$\bullet$};
            \node at (8.7,-6) {$10$};
            \node at (11,-6) {$\bullet$};
            \node at (11.3,-6) {$00$};
            \node at (10,-7) {$\bullet$};
            \node at (10,-7.3) {$00$};

            \draw  (9,-6) -- (11,-5) -- (11,-6) -- (9,-5) -- (9,-6) -- (10,-7) -- (11,-6);

            \node at (12,-5) {$\bullet$};
            \node at (12,-4.7) {$11$};
            \node at (14,-5) {$\bullet$};
            \node at (14,-4.7) {$10$};
            \node at (13,-6) {$\bullet$};
            \node at (13.3,-6) {$10$};
            \node at (13,-7) {$\bullet$};
            \node at (13,-7.3) {$00$};

            \draw (12,-5) -- (13,-6) -- (14,-5) -- (13,-6) -- (13,-7);

            \node at (2,-11) {$Q_{5}$};

            \node at (5,-9) {$\bullet$};
            \node at (5,-8.7){$111$};
            \node at (4,-10) {$\bullet$};
            \node at (3.7,-10){$110$};
            \node at (4,-11) {$\bullet$};
            \node at (3.7,-11){$100$};
            \node at (6,-10) {$\bullet$};
            \node at (6.3,-10){$100$};
            \node at (6,-11) {$\bullet$};
            \node at (6.3,-11){$000$};
            \node at (5,-12) {$\bullet$};
            \node at (5,-12.3) {$000$};

            \draw (5,-9) -- (4,-10) -- (4,-11) -- (6,-10) -- (4,-11) -- (4,-10) -- (6,-11) -- (6,-10) -- (5,-9) -- (4,-10) -- (4,-11) -- (5,-12) -- (6,-11) -- (4,-10);

            \node at (8,-11) {$\bullet$};
            \node at (8.3,-11) {$100$};
            \node at (8,-12) {$\bullet$};
            \node at (8,-12.3) {$000$};
            \node at (7,-10) {$\bullet$};
            \node at (7.3,-10) {$110$};
            \node at (9,-10) {$\bullet$};
            \node at (9.3,-10) {$100$};
            \node at (8,-9) {$\bullet$};
            \node at (8,-8.7) {$111$};

            \draw (8,-9) -- (7,-10) -- (8,-11) -- (9,-10) -- (8,-9) -- (9,-10) -- (8,-11) -- (8,-12);

        \end{tikzpicture}
        \caption{Models witnessing failure of $3$-uniformity}
        \label{fig:2bisimilar models not bisimilar}
    \end{figure}

\subsection{A Non-Uniform Locally Tabular Logic}\label{section: non-uniform locally tabular logic}

The analysis of $3$-uniformity sheds some light on the possibilities for logics to be $n$-uniform but not $(n+1)$-uniform. Namely, note that $\mathsf{Log}(Q_{4})$ will be an extension of $\mathsf{wPL}$, and it will be $4$-uniform but not $3$-uniform by what was just shown. A more general fact holds which will allow us to provide an example of a locally tabular logic which is not uniformly locally tabular.

\begin{definition}
    Let $n>2$ and $k+1<n$. Let $M_{n}^{k}$, $N_{n}^{k}$ for $k\leq n$ be the models depicted in Figure \ref{fig:modelsofcascades}. Formally, we define:
    \begin{enumerate}
        \item $M_{n}^{k}$ has $k+1$ many layers of two elements; if $L_{j}(M_{n}^{k})$ is the $j$-th layer for $0\leq j\leq k$, then the left element in $L_{j}(M_{n}^{k})$ is labelled with $1^{n-j}0^{j}$, and the right element with $1^{n-(j+1)}0^{j+1}$. There is also a root element labelled with $1^{n-k-1}0^{k+1}$.
        \item $N_{n}^{k}$ has likewise $k+1$ many layers of two elements; if $L_{j}(N_{n}^{k})$ is the $j$-th layer for $0\leq j\leq k$, then the left element in $L_{j}(N_{n}^{k})$ is labelled with $1^{n-j}0^{j}$, and the right element with $1^{n-(j+1)}0^{j+1}$. There is also a root element labelled with $1^{n-k-2}0^{k+2}$.
    \end{enumerate}
    Thus the only difference between $M_{n}^{k}$ and $N_{n}^{k}$ is in the label of the root. The order in both cases is given by stacking the layers on top of each other.
\end{definition}

Note that these models are always $n$-generated. Indeed, the key difference between these models is that the root of $N_{n}^{k}$ goes down in color, relative to its right immediate successor, whilst in $M_{n}^{k}$ the color is the same. Figure \ref{fig:modelsm1andn0} shows the cases for $k=1$, for any $n\geq 2$

\begin{figure}[h]
    \centering
\begin{tikzpicture}
    \node at (0,0) {$\bullet$};
    \node at (0,0.3) {$1^{n}$};
    \node at (2,0) {$\bullet$};
    \node at (2,0.3) {$1^{n-1}0$};
    \node at (0,-1) {$\bullet$};
    \node at (0,-1.3) {$1^{n-1}0$};
    \node at (2,-1) {$\bullet$};
    \node at (2,-1.3) {$1^{n-2}00$};
    \node at (1,-2) {$\bullet$};
    \node at (1,-2.3) {$1^{n-2}00$};
    \node at (1,-3) {$M_{n}^{1}$};

    \draw (0,0) -- (2,-1) -- (2,0) -- (0,-1) -- (0,0) -- (0,-1) -- (1,-2) -- (2,-1);

    \node at (4,0) {$\bullet$};
    \node at (4,0.3) {$1^{n}$};
    \node at (6,0) {$\bullet$};
    \node at (6,0.3) {$1^{n-1}0$};
    \node at (5,-1) {$\bullet$};
    \node at (5,-1.3) {$1^{n-2}00$};
    \node at (5,-3) {$N_{n}^{0}$};

    \draw (4,0) -- (5,-1) -- (6,0);
\end{tikzpicture}    \caption{Models for $\mathsf{wPL}$}
\label{fig:modelsm1andn0}
\end{figure}

The following lemma provides the relationship between such models:

\begin{figure}
    \centering
\begin{tikzpicture}
    \node at (0,0) {$\bullet$};
    \node at (0,0.3) {$1^{n}$};
    \node at (0,-1) {$\bullet$};
    \node at (0,-1.3) {$1^{n-1}0$};
    \node at (2,0) {$\bullet$};
    \node at (2.3,0.3) {$1^{n-1}0$};
    \node at (2,-1) {$\bullet$};
    \node at (2.3,-1.3) {$1^{n-2}00$};
    \node at (1,-2) {$\vdots$};
    \node at (0,-4) {$\bullet$};
    \node at (0,-3.7) {$1^{n-k}0^{k}$};
    \node at (2,-4) {$\bullet$};
    \node at (2,-3.7) {$1^{n-k-1}0^{k+1}$};
    \node at (0,-5) {$\bullet$};
    \node at (-0.3,-5.3) {$1^{n-k-1}0^{k+1}$};
    \node at (2,-5) {$\bullet$};
    \node at (2,-5.3) {$1^{n-k-2}0^{k+2}$};
    \node at (1,-6) {$\bullet$};
    \node at (1.2,-6.3) {$1^{n-k-2}0^{k+2}$};
    \node at (1,-7.5) {$M_{n}^{k+1}$};

    \draw (1,-6) -- (0,-5) -- (0,-4) -- (2,-5) -- (2,-4) -- (0,-5) -- (1,-6) -- (2,-5);
    
    \draw (0,0) -- (2,-1) -- (2,0) -- (0,-1) -- (0,0);

    \node at (4,0) {$\bullet$};
    \node at (4,0.3) {$1^{n}$};
    \node at (4,-1) {$\bullet$};
    \node at (4,-1.3) {$1^{n-1}0$};
    \node at (6,0) {$\bullet$};
    \node at (6.3,0.3) {$1^{n-1}0$};
    \node at (6,-1) {$\bullet$};
    \node at (6.3,-1.3) {$1^{n-2}00$};
    \node at (5,-2) {$\vdots$};
    \node at (4,-4) {$\bullet$};
    \node at (4,-3.7) {$1^{n-k+1}0^{k-1}$};
    \node at (6,-4) {$\bullet$};
    \node at (6,-3.7) {$1^{n-k}0^{k}$};
    \node at (4,-5) {$\bullet$};
    \node at (4,-5.3) {$1^{n-k}0^{k}$};
    \node at (6,-5) {$\bullet$};
    \node at (6,-5.3) {$1^{n-k-1}0^{k+1}$};
    \node at (5,-6) {$\bullet$};
    \node at (5.2,-6.3) {$1^{n-k-2}0^{k+2}$};
    \node at (5.2,-7.5) {$N_{n}^{k}$};

    \draw (5,-6) -- (4,-5) -- (4,-4) --(6,-5) -- (6,-4) -- (4,-5) -- (5,-6) -- (6,-5);
    \draw (4,0) -- (6,-1) -- (6,0) -- (4,-1) -- (4,0);
\end{tikzpicture}

\caption{General models $M_{n}^{k+1}$ and $N_{n}^{k}$}    \label{fig:modelsofcascades}
\end{figure}

\begin{lemma}\label{N-bisimulation lemma}
    Let $n>2$ be arbitrary, and $k+1<n$. Then $M_{n}^{k}$ is $k$-bisimilar, but not $(k+1)$-bisimilar with $N_{n}^{k-1}$.
\end{lemma}
\begin{proof}
    This is shown by induction. Note that for $k=1$, this is easy to verify, since $N_{n}^{1}$ will be of depth $2$: looking at Figure \ref{fig:modelsm1andn0}, we have that the two roots have the same color, and they see points with the same colors, so they will be $1$-bisimilar. However, assume that they were $2$-bisimilar; let $x$ be the point in $M_{n}^{1}$ with the label $1^{n-1}0$ which is the successor of the root. Then by $2$-bisimilarity, there is some point in $N_{n}^{0}$ which is $1$-bisimilar with it; but the only point which could be is the maximal point, which is a contradiction.

    Now, assume that this holds for $k$. Consider $M_{n}^{k+1}$ and $N_{n}^{k}$. Note that they have roots of the same color, labelled with $1^{n-k-2}0^{k+2}$, and indeed, if $z$ is the immediate successor of $M_{n}^{k+1}$ with the same color, ${\uparrow}z$ is isomorphic to $N_{n}^{k}$. Thus if Player I picks the root of $N_{n}^{k}$ or any of its successors, Player II has an obvious response by choosing the element inside of ${\uparrow}z$. Thus the only way that Player I can hope to win the game is to pick either the root of $M_{n}^{k}$, or the left successor of the root. If Player I picks the left successor of the root, say $y_{0}$, note that ${\uparrow}y_{0}\cong M_{n}^{k}$; on the other hand, the right successor of the root of $N_{n}^{k}$ is isomorphic to $N_{n}^{k-1}$, so Player II wins by induction hypothesis. Finally, if Player I picks the root of $M_{n}^{k+1}$, Player II is forced to show that this is $k$-bisimilar to the root of $N_{n}^{k}$; by the arguments just shown Player I is in turn forced to play the root of $M_{n}^{k+1}$ for all $k$ many steps of the game, and loses when Player II finally picks the root, given the two roots have the same color. This shows that $M_{n}^{k+1}$ is $k$-bisimilar to $N_{n}^{k}$ as desired.

    Now assume that $M_{n}^{k+1}$ is $(k+1)$-bisimilar to $N_{n}^{k}$. Then Player I can pick the left successor of the root of $M_{n}^{k+1}$, $y_{0}$, labelled with $1^{n-k-1}0^{k+1}$, and Player II is forced to respond with the right successor of $N_{n}^{k}$, $y_{1}$; as noted ${\uparrow}y_{0}\cong M_{n}^{k}$ and ${\uparrow}y_{1}\cong N_{n}^{k-1}$. Since by induction hypothesis these models are not $k$-bisimilar, Player II then loses, meaning that $M_{n}^{k+1}$ and $N_{n}^{k}$ are not $(k+1)$-bisimilar.
\end{proof}

Hence we have the following:

\begin{theorem}\label{thm: loc tabular but not uniformly so}
    The logic $\mathsf{wPL}$ is locally tabular but not uniformly locally tabular.
\end{theorem}
\begin{proof}
    All the models $M_{n}^{k}$ for each $n,k$ are based on frames which validate $\mathsf{wPL}$; hence by Corollary \ref{Saturated classes yielding n-uniformity}, $\mathsf{wPL}$ is not uniformly locally tabular, since we notice that $M_{n}^{k+1}\leq_{k+1}N_{n}^{k}$ but not $M_{n}^{k+1}\leq_{k+2}N_{n}^{k}$. On the other hand, it is locally tabular by Theorem \ref{Local tabularity of wPL}.
\end{proof}

\section{Pre-uniform local tabularity}\label{sec: preuniform local tabularity}

As noted before the problem of working with locally tabular intermediate logics presents a challenge since, as shown by Mardaev \cite{Mardaev1984}, there are uncountably many such pre-locally tabular logics. By contrast, there are exactly three pretabular logics \cite{ESAKIA1977}, and every non-tabular logic is contained in one of them, and all of them are decidable. Together, this makes the problem of deciding, given an axiom $\phi$, whether $\mathsf{IPC}\oplus \phi$ is tabular or not decidable. We thus have the following interesting problem:

\begin{problem}\label{big problem}
    Is uniform local tabularity decidable over $\mathsf{IPC}$?
\end{problem}

Intrinsic in this is the idea of ``pre-local tabularity''.  The idea of ``pre-$P$''-logics, for $P$ a property, has been extensively used by several authors (e.g., \cite{ESAKIA1977,Segerberg1971-SEGAEI,Maksimova1975}) to characterize certain logical properties. We can make such a definition also in our setting:

\begin{definition}
    Let $L$ be a superintuitionistic logic. We say that $L$ is \textit{pre-uniformly locally tabular} if $L$ is not uniformly locally tabular but all of its extensions are so.
\end{definition}

We will now characterize pre-uniform local tabularity amongst extensions of $\mathsf{KG}$. For that, first recall the following well-known formula,
\begin{equation*}
    \mathsf{bw}_{2}\coloneqq (p_{0}\rightarrow (p_{1}\vee p_{2}))\vee (p_{1}\rightarrow (p_{0}\vee p_{2}))\vee (p_{2}\rightarrow (p_{0}\vee p_{1})),
\end{equation*}
the \emph{bounded width 2} formula \cite[Proposition 2.39]{Chagrov1997-cr}; a finite frame validates this axiom if and only if they have width at most $2$, i.e., for each point $x$ in such a frame, ${\uparrow}x$ contains antichains of size at most $2$. Recall also the formula:
\begin{equation*}
    \neg p\vee \neg\neg p,
\end{equation*}
the weak law of excluded middle; it is also well-known \cite[Proposition 2.37]{Chagrov1997-cr} that frames for such a formula are exactly those containing a greatest point above each point.

\begin{definition}
    The logic $\mathsf{Box}$ is defined to be the extension of $\mathsf{wPL}$ by the  $\mathsf{Box}=\mathsf{wPL}\oplus \mathsf{bw}_{2}\oplus (\neg p\vee \neg\neg p)$.
\end{definition}

We will see that rooted frames for $\mathsf{Box}$ look very much like those from Figure \ref{fig:modelsofcascades}, with a point in the maximum --- if one thinks of four adjacent points there, these look like stacks of boxes.

\begin{remark}
    Note that $\mathsf{Box}$, as an extension of $\mathsf{wPL}$ is locally tabular. For that reason, we will assume throughout that all logics $L$ extending it are locally tabular, and consequently have the finite model property.
\end{remark}

We will now show that $\mathsf{Box}$ is pre-uniformly locally tabular. For that purpose, let $S_{n}$ be the result of adding a maximal point to the underlying frame of $N_{n}^{n-2}$ (see Figure \ref{fig:modelsofcascades}). We denote by $\mathsf{Lev}_{j}(S_{k})=\{x\in S_{k}:d_{p}(x)=j\}$, the \emph{j-th level}. Then first we show:

\begin{proposition}\label{prop: Generation of Box logic}
    We have that $\mathsf{Box}=\mathsf{wPL}\oplus \mathsf{bw}_{2}\oplus (\neg p\vee \neg \neg p)=\mathrm{Log}(\{S_{n} : n\in \omega\})$.
\end{proposition}
\begin{proof}
    It is clear that each $S_{n}$ validates all the axioms of $\mathsf{Box}$, since they have a maximal point, have width at most $2$, and are frames for $\mathsf{wPL}$ as noted in Theorem \ref{thm: loc tabular but not uniformly so}. Conversely, by the semantic characterizations given of these formulas, any finite rooted poset $P$ validating $\mathsf{Box}$ will be a Boolean sum (as it validates $\mathsf{wPL} $), it will have width at most $2$ by $\mathsf{bw}_{2}$, and a maximal element by the axiom $\neg p\vee \neg\neg p$. Let $k$ be the depth of the root of $P$. Then for any $j\leq k$, either there is a unique point $x$ such that $d_{p}(x)=j$, or there are two. This defines a p-morphism from $S_{k}$ onto $P$. Thus, every finite rooted poset of $\mathsf{Box}$ can be obtained by p-morphisms from one of the form $S_{n}$, showing the result.
\end{proof}

Recall the notion of stack depth (Definition \ref{def: stack depth}); for frames with bounded width $2$, having stack depth $k$ amounts to having a bound of $k$ on the number of consecutive levels with two points. Thus we obtain:

\begin{proposition}
    Let $\mathsf{Box}\subsetneq L'$ be a logic of infinite depth. Then there is some $n$ such that posets in $L'$ have stack depth $n$.
\end{proposition} 
\begin{proof}
    If we assume that there is no bound on the stack depth of posets for $L'$, then for each $n$, let $P$ be a frame for $L'$ containing a stack of at least $n+2$ many elements. Let $p$ be one of the minimal elements of the stack, and $m$ be an element situated at at $n$ many levels above $p$; then ${\uparrow}p$ is a frame for $L'$. Define a p-morphism onto $S_{n}$ by sending $m$, the other element in its level, and all their successors to the top element in $S_{n}$, and all elements in the stack mapped through an isomorphism to $S_{n}$. Thus $S_{n}\Vdash L'$. Since this holds for any $n$, by  Proposition \ref{prop: Generation of Box logic} we have that $L'=\mathsf{Box}$.
\end{proof}

We thus obtain that the extensions of $\mathsf{Box}$ are controlled by the depth, which can be infinite or finite, and the stack-depth, which can be finite or infinite. We show the following:

\vspace{2mm}

\begin{proposition}\label{Bound on stack depth lemma}
    Assume that $L\supseteq \mathsf{Box}$ is a logic of infinite depth such that all its elements have at most stack-depth $k\geq 1$. Then $L$ is $(k+2)$-uniformly tabular. Consequently, $\mathsf{Box}$ is pre-uniformly locally tabular.
\end{proposition}
\begin{proof}
The fact that $\mathsf{Box}$ is not uniformly locally tabular follows by a simple modification of Theorem \ref{thm: loc tabular but not uniformly so}, adding a point on top to the models $M_{n}^{k}$ with a distinguished color.

Assume that $L\supseteq \mathsf{Box}$ is a logic such that its models have stack depth at most $k$. Then we show the following modification of Lemma \ref{lem: 1-bisim}: if $(P,x)$ is a Boolean sum of stack depth at most $k$ and $(P,x)\sim_{k}(P_{y},y)$ then $x=y$. This proceeds by a double induction on $k$ and $m$, where $m$ is the number of stacks above $x$ of size at least $k$. For $k=1$ and arbitrary $m$, this was provided by Lemma \ref{lem: 1-bisim}. Assume that it holds for $k$, and we first show it for $k+1$ and $m=1$; In this case, note that since $(P,x)\sim_{k+1}(P_{y},y)$, we can pick an immediate successor of $x$, $x'$, and by $k$-bisimilarity, $(P_{x'},x')\sim_{k}(P_{y'},y')$ for some sucessor of $y$; then the stack above $P_{x'}$ must be of size at most $k$, so by induction hypothesis, $x'=y'$ must hold. Since this holds for any immediate successor of $x$, we conclude that $x=y$. Assuming this for arbitrary $m$, we then show it for $m+1$, by noting that when we move to an immediate successor, the immediate successor will see one fewer $k+1$-stack, and so the induction hypothesis can be carried through.

From this fact, we proceed as in Proposition \ref{Characterization of 2-uniformity}: assuming we have $(P,x),(Q,y)$, $n$-generated models based on frames for $L$, where $L$ has a bound on the stack-depth of $k$, such that $(P,x)\sim_{k+1}(Q,y)$. Then we have an $S_{k+1}$-bisimulation, defined by a chain $S_{k+1}\subseteq S_{k}\subseteq\dots \subseteq S_{0}$; then the relation $S_{k}$ defines an isomorphism, and so $(P,x)$ and $(Q,y)$ are bisimilar. From this it follows by the same arguments as in Proposition \ref{Characterization of 2-uniformity} that $L$ is $k+2$-uniformly locally tabular.
\end{proof}

\subsection{Pre-uniform local tabularity in KG}

We will provide a solution to Problem \ref{big problem} in the case of extensions of the system $\mathsf{KG}$ of Kuznetsov-Gerciu \cite{kuznetsovgerciu} (for an overview of the system see \cite[Chapter 4]{bezhanishviliphdthesis}). Such a system has a very broad range of extensions, and a rich semantics, which has made it useful for several purposes in the study of superintuitionistic logic (e.g. \cite{Bezhanishvili2025}). One of its key extensions is the logic $\mathsf{RN}=\mathrm{Log}(\mathsf{RN})$ --- the logic of the frame of the Rieger-Nishimura ladder (see Figure \ref{fig:riegernishimura}), which is dual to the \textit{Rieger-Nishimura lattice}, the free $1$-generated Heyting algebra.

\begin{figure}[h]
    \centering
\begin{tikzpicture}
    \node at (0,0) {$\bullet$};
    \node at (0,0.3) {$1$};
    \node at (1,0) {$\bullet$};
    \node at (1,0.3) {$0$};
    \node at (0,-1) {$\bullet$};
    \node at (-0.3,-1) {$P_{1}$};
    \node at (1,-1) {$\bullet$};
    \node at (1.3,-1) {$P_{2}$};
    \node at (0,-2) {$\bullet$};
    \node at (-0.3,-2) {$P_{3}$};
    \node at (1,-2) {$\bullet$};
    \node at (1.3,-2) {$P_{4}$};
    \node at (0,-3) {$\bullet$};
    \node at (-0.3,-3) {$P_{5}$};
    \node at (1,-3) {$\bullet$};
    \node at (1.3,-3) {$P_{6}$};
    \node at (0,-4) {$\bullet$};
    \node at (-0.3,-4) {$P_{7}$};
    \node at (1,-4) {$\bullet$};
    \node at (1.3,-4) {$P_{8}$};
    \node at (0.5,-5) {$\vdots$};
    \node at (0.5,-6) {$\bullet$};
    \node at (0.5,-6.3) {$\infty$};

    \draw (0,0) -- (0,-1) -- (0,0) -- (1,-1) -- (1,0) -- (1,-1) -- (1,-2) -- (0,-1) -- (0,-2) -- (1,0);

    \draw (0,-2) -- (0,-3) -- (1,-1) -- (1,-2) -- (1,-3) -- (0,-2);

    \draw (0,-3) -- (0,-4) -- (1,-2) -- (1,-3) -- (1,-4) -- (0,-3);
\end{tikzpicture}
    
    \caption{Rieger-Nishimura Ladder $\mathsf{RN}$}
    \label{fig:riegernishimura}
\end{figure}

Recall from \cite{ghilardifreeheyting,Butz1998} (see also \cite{Bezhanishvili2011} and \cite{colimitsofHeytingalgebras}) that $\mathsf{RN}$, the Rieger-Nishimura ladder, can be constructed as the inverse limit of a sequence:

\begin{equation*}
    P_{1}\xleftarrow{r_{1}} P_{2}\xleftarrow{r_{2}} P_{3}\xleftarrow{r_{3}}\dots
\end{equation*}
where $r_{i}$ is a map which is a p-morphism with respect to ``sufficiently many'' parts of $P_{3}$. In fact these admit an explicit description:
\begin{enumerate}
    \item The posets $P_{i}$ are precisely the upsets as labelled in Figure \ref{fig:riegernishimura}.
    \item The map $r_{2n+1}:P_{2n+2}\to P_{2n+1}$ sends the $P_{2n+2}$-root point, and the $P_{2n}$-root point, to $P_{2n+1}$, and leaves the rest untouched;
    \item The map $r_{2n+2}$ sends the $P_{2n+3}$ and the $P_{2n+1}$-root points to $P_{2n+2}$,  and leaves everything else untouched. 
\end{enumerate}

For our purposes we will need the following lemma, shown essentially in \cite[Chapter 4]{bezhanishviliphdthesis}:

\begin{lemma}\label{lem: obtaining the boxes out of p-morphisms}
    For each $n$, the poset $S_{n}$ is an upset of a p-morphic image of $P_{2n+1}\oplus 1$.
\end{lemma}
\begin{proof}
    One proceeds inductively as follows: first identify the point labelled $1$ in Figure \ref{fig:riegernishimura} with the top element. This produces one top level containing the maximum, followed by a level containing $P_{1}$ and $0$, seen by $P_{3}$. Then to obtain the next layer, identify $P_{2}$ with its now unique successor, and we obtain a new layer, consisting of $P_{3}$ and $P_{4}$, seen by $P_{6}$. We thus proceed in this way recursively, obtaining all boxes as desired.
\end{proof}

To characterize uniform local tabularity, we will rely on the fact that already local tabularity above $\mathsf{KG}$ is characterized. In principle this is not necessary, as uniform local tabularity is potentially simpler; however in this simple case, it will be convenient to rely on the following, shown in \cite[Theorem 4.6.5]{bezhanishviliphdthesis}:

\begin{theorem}\label{thm: criterion for prelocal tabularity}
    A logic $\mathsf{KG}\subseteq L$ is not locally tabular if and only if $L\subseteq \mathsf{RN.KC}=\mathsf{RN}\oplus (\neg p\vee \neg \neg p)$.
\end{theorem}

We will also need one more technical notion:

\begin{definition}
    We say that a poset $P$ has \emph{degree of uniformity $k$} if whenever $(P,x)$ is an $n$-generated model, and $(Q,y)$ is another $n$-generated model such that $S_{k}\subseteq P\times Q$, is a $k$-bisimulation between $x$ and $y$, then $S_{k}$ is an isomorphism.
\end{definition}

\begin{remark}
    The reader should notice that for a fixed poset $P$, a degree of uniformity need not be minimal --- for the purposes where we will use it, the notion of a minimal degree of uniformity will not be needed.
\end{remark}

It is well-known that any finite poset has a degree of uniformity (see e.g. \cite[Proposition 4.5]{Ghilardi2002}). Recall that given rooted posets $P,Q$, we write $P\oplus Q$ for the \emph{vertical sum} (see \cite[Def.~4.1.13]{bezhanishviliphdthesis}) the poset based on the disjoint union of $P,Q$, where for each $x\in P$ and $y\in Q$, $x\leq y$. We denote the repeated vertical sum as $P_{1}\oplus\dots \oplus P_{n}$.

\begin{lemma}\label{lem: uniform local tabularity from sums}
    If $P_{1},\dots,P_{n}$ is a sequence of posets with degree of uniformity $k_{i}$, then $P_{1}\oplus\dots\oplus P_{n}$ has degree of uniformity $k'=\max(\{k_{i} : i\leq n\})+1$. 
\end{lemma}
\begin{proof}
    Write $P=P_{1}\oplus\dots\oplus P_{n}$. Let $(Q,y)$ be an arbitrary $n$-generated model, and assume that $(P,x)\sim_{k'}(Q,y)$. Let $x_{n}$ be the root of $P_{n}$; by $k'$-bisimilarity, we have that there is some $y_{n}\in Q$, such that $(P_{x_{n}},x_{n})\sim_{k'-1}(Q_{y_{n}},y_{n})$; since $k_{n}\leq k'-1$, this relation defines an isomorphism. Note that any predecessor of $x_{n}$ or $y_{n}$ must have a properly smaller color. Let $Q_{n}={\uparrow}y_{n}$. By restricting the relation $S_{k'}$ to $P\setminus (P_{n}\setminus \{x_{n}\})$ and $Q\setminus (Q_{n}\setminus\{y_{n}\})$, we claim that we still have a $k'$-bisimulation: for each $x'\geq x$, such that $x'\in P\setminus (P_{n}\setminus \{x_{n}\})$, there must be $y'\geq y$ which is $k'-1$-bisimilar, and this must lie also in $Q\setminus (Q_{n}\setminus\{y_{n}\})$. Thus we can proceed by induction, showing that each such restriction of $S_{k}$ is an isomorphism.
\end{proof}

  \begin{figure}[h]
    \centering
\begin{tikzpicture}
\filldraw[color=orange!50, fill=orange!5,thick] (0,1) circle (4cm);

\filldraw[color=blue!50, fill=blue!5, thick] (-1,1) ellipse (1.6 and 3.5);

    \draw (0,0) circle (5cm);
    
    \draw (0,2) circle (3cm);
    \draw (0,2.75) circle (2.25cm);
    \draw (-1,3.25) circle (1cm);

    \node at (-7,-4) {$\mathbf{LTab}$};

\draw[dotted] (-6.3,-4) -- (0,-4);
    
    \node at (0,-4.5) {$\mathsf{wPL}$};
    \node at (0,-3.15) {$\mathsf{Box}$};
    \node at (0,-2) {$\vdots$};

    \node at (-1,-2) {$\mathsf{Log}(P_{4})$};
    
    \node at (-7,-0.5) {$3$-$\mathbf{Unif}$};

    \draw[dotted] (-6.3,-0.5) -- (1,-0.5);

    \node at (-7,1.5) {$2$-$\mathbf{Unif}$};
    \node at (-1,0.3) {$\mathsf{Log}(P_{2})$};
    \node at (1,0) {$\mathsf{2Uni}$};
    \node at (1,2.5) {$\mathsf{LC}$};

    \draw[dotted] (-6.3,1.5) -- (1.5,1.5);

    \node at (-7,3.25) {$1$-$\mathbf{Unif}$};
    \node at (-1,3.5) {$\mathsf{CPC}$};

    \draw[dotted] (-6.3,3.25) -- (-0.75,3.25);

    \node at (7,0.25) {$\mathbf{Tab}$};

    \draw[dotted] (6.3,0.25) -- (0,0.25);

    \node at (7,-2.25) {$\mathbf{ULTab}$};

    \draw[dotted] (6.3,-2.25) -- (0.25,-2.25);

\end{tikzpicture}    \caption{Hierarchy of Locally Tabular Superintuitionistic Logics; blue oval covers the tabular logics.}
    \label{fig:hierarchyofsuperintuitionisticlogics}
\end{figure}

We can thus conclude:

\begin{theorem}\label{thm: characterization of uniform local tabularity above KG}
If $L\supseteq \mathsf{KG}$ is any logic, then $L$ is not pre-uniformly locally tabular if and only if $L\subseteq\mathsf{Box}$.
\end{theorem}
\begin{proof}
    If $\mathsf{KG}\subseteq L$, and $L$ is not locally tabular, we thus have by Theorem \ref{thm: criterion for prelocal tabularity}  that $L\subseteq \mathsf{RN}\subseteq \mathsf{Box}$. Hence assume that $L$ is locally tabular but not uniformly locally tabular. We can assume by Theorem \ref{prop: Inclusion of classes of logics} that $L$ is not of finite depth. Since $L$ is locally tabular, there must be an $n$ which is maximal such that ${\uparrow}P_{n}'$ is a frame for the logic. Hence by \cite[Corollary 4.3.9]{bezhanishviliphdthesis}, the finitely generated frames for $L$ must be finite vertical sums of p-morphic images of upsets of $P_{n}$ for $k\leq n$ (not necessarily rooted).

    Now as noted in Lemma \ref{lem: obtaining the boxes out of p-morphisms}, for large enough $P_{n}$ we can obtain $S_{n}$, and for large enough vertical sums, we can obtain all the $S_{n}$. Suppose that there is some $n$ such that we cannot obtain $S_{n}$ as a frame for $L$. Then there must be a bound on the size of the posets $Q$ with a root and a maximal element which appear in vertical sums in $L$ (otherwise using the technique of Lemma \ref{lem: obtaining the boxes out of p-morphisms} we could obtain arbitrarily large $S_{n}$). Using Lemma \ref{lem: uniform local tabularity from sums} we can conclude that posets for $L$ will have a maximal degree of uniformity $k$ for some $k$ --- which implies (and is stronger than the fact that) whenever $(P,x)\leq_{k}(Q,y)$ for $(P,x)$ and $(Q,y)$ two models based on frames for $L$, then $(P,x)\leq_{\infty}(Q,y)$. By Proposition \ref{Saturated classes yielding n-uniformity} we obtain that $L$ is uniformly locally tabular --- a contradiction.

    So we must assume that $L$ contains arbitrarily large vertical sums of $P_{n}$ as frames. Thus as noted, we can obtain $S_{n}$ for each $n$. By Proposition \ref{prop: Generation of Box logic} we have that $L\subseteq \mathsf{Box}$ as desired.
\end{proof}

\begin{corollary}
    Above $\mathsf{KG}$, uniform local tabularity is decidable: if $L=\mathsf{KG}\oplus \phi$, it is decidable whether $L$ is uniformly locally tabular.
\end{corollary}
\begin{proof}
    Since the logic $\mathsf{Box}$ is finitely axiomatized and locally finite, it is decidable. Given $L$ a decidable logic, one simply checks whether $L\subseteq \mathsf{Box}$, by checking whether $\mathsf{Box}$ contains the axiom $\phi$, which is enough by Theorem \ref{thm: characterization of uniform local tabularity above KG}. 
\end{proof}

\section{Summary}

As a conclusion to the paper, we note a visual summary of our results in Figure \ref{fig:hierarchyofsuperintuitionisticlogics}.

\section{Acknowledgements}

The author would like to thank Nick Bezhanishvili, Luca Carai and Tommaso Moraschini for valuable discussions on the results of this paper, as well as the anonymous referees for valuable suggestions that greatly improved the presentation of the present paper, and for pointing some imprecisions which have been corrected in the present version.


\newpage

\end{document}